\documentclass[12pt]{amsart}
\usepackage{amsfonts}
\usepackage{amsmath}
\usepackage{amssymb}

\newtheorem{theorem}{Theorem}[section]

\newtheorem{conjecture}[theorem]{Conjecture}
\newtheorem{corollary}[theorem]{Corollary}

\newtheorem{definition}[theorem]{Definition}

\newtheorem{lemma}[theorem]{Lemma}

\newtheorem{proposition}[theorem]{Proposition}
\newtheorem{remark}[theorem]{Remark}

\numberwithin{equation}{section}

\begin{document}
\title{Supremum of Perelman's  entropy  and K\"ahler-Ricci flow on a Fano manifold}

\author   {Gang Tian}
\author{Shijin Zhang}
\author{Zhenlei  $\text{Zhang}^*$}
\author { Xiaohua $\text{Zhu}^{**}$}

\thanks {*  Partially supported by  a grant  of BMCE 11224010007 in China}
\thanks {**  Partially supported by the NSFC Grant 10990013}

\address {Gang Tian\\ School of Mathematical Sciences and BICMR, Peking
University, Beijing, 100871, China\\ and Department of Mathematics, Princeton
University,  New Jersey,  NJ 02139, USA\\
 tian@math.mit.edu}
\address{Shijin Zhang, School of Mathematical Sciences and BICMR, Peking University, Beijing, 100871, P.R.China, shijin\_zhang@yahoo.com}
\address{Zhenlei   Zhang\\ Department of Mathematics,  Beijing Capital Normal University, Beijing, China\\     zhleigo@yahoo.com.cn}
\address{ Xiaohua Zhu\\School of Mathematical Sciences and BICMR, Peking University,
Beijing, 100871, China\\
 xhzhu@math.pku.edu.cn}

\maketitle

\begin{abstract}
In this paper, we extend  the method in [TZhu5] to study the energy level $L(\cdot)$ of Perelman's entropy
$\lambda(\cdot)$  for  K\"ahler-Ricci flow on a Fano manifold.  Consequently,
we  first compute the supremum of 
$\lambda(\cdot)$ in  K\"ahler class $2\pi c_1(M)$  under an assumption  that the modified Mabuchi's K-energy
$\mu(\cdot)$ defined in [TZhu2] is bounded from below.  Secondly,  we give an alternative proof to the main theorem about the convergence of
 K\"ahler-Ricci flow in [TZhu3].
\end{abstract}

\vskip6mm
\section*{Introduction}

In this paper, we extend  the method in [TZhu5] to study the energy level $L(\cdot)$ of Perelman's entropy
$\lambda(\cdot)$  for  K\"ahler-Ricci flow on an $n$-dimensional  compact K\"{a}hler manifold $(M,J)$
with positive first Chern class $c_1(M)>0$ (namely  called a Fano
manifold).  We will show  that  $L(\cdot)$ is independent of choice of initial K\"ahler metrics in $2\pi c_1(M)$ under an assumption that
the modified Mabuchi's K-energy  $\mu(\cdot)$ is bounded from below (cf. Proposition \ref{prop-2}  in Section 3).  The modified Mabuchi's K-energy  $\mu(\cdot)$ is a generalization of
Mabuchi's K-energy. It was showed in [TZhu2] that  $\mu(\cdot)$  is
bounded from below if   $M$ admits a K\"{a}hler-Ricci soliton.

 As an application of  Proposition \ref{prop-2}, we first
 compute the supremum of  Perelman's entropy $\lambda(\cdot)$ in K\"ahler class $2\pi c_{1}(M)$ [Pe].  More precisely, we  prove that

\begin{theorem}\label{theorem-1}
Suppose that the modified Mabuchi's K-energy is bounded from below.
Then
\begin{equation}
\sup\{\lambda({g'})|~g'\in
\mathcal{K}_{X}\}=(2\pi)^{-n}[nV-N_X(c_1(M))].
\end{equation}
\end{theorem}

Here the quantity $N_X(c_1(M))$  is a  nonnegative invariance in
$\mathcal{K}_{X}$  and it is zero iff  the  Futaki-invariant  vanishes [Fu].    We denote $\mathcal{K}_{X}$  to be a  class of $K_X$-invariant  K\"{a}hler metrics in $2\pi c_{1}(M)$,
where  $K_X$  is an  one-parameter  compact subgroup  of holomorphisms transformation group   generated by  an extremal  holomorphic vector field $X$  for
K\"ahler-Ricci solitons on $M$ [TZhu2]. We note that we  do not  need to assume  an existence of K\"ahler-Ricci solitons in Theorem \ref{theorem-1}.
 In fact,  if we assume the existence of K\"ahler-Ricci solitons   then we can  use a more direct
way to prove Theorem \ref{theorem-1}  and   that  the supremum of $\lambda(\cdot)$ can be achieved
in $\mathcal{K}_{X}$ (cf. Section 1).   It seems that the supremum of
$\lambda(\cdot)$ can  be achieved  in the total space of K\"ahler
potentials in $2\pi c_1(M)$ if $M$ admits a K\"{a}hler-Ricci
soliton. In a  special case of small neighborhood of a
K\"{a}hler-Ricci soliton the positivity has been verified  by
computing the second variation of $\lambda(\cdot)$ in [TZhu4].

 As  another application of  Proposition \ref{prop-2},  we prove the  following  convergence  result  for K\"ahler-Ricci flow.

\begin{theorem}\label{theorem-2} Let $(M,J)$ be a compact K\"ahler manifold
which admits a  K\"ahler-Ricci soliton $(g_{KS},X)$.   Then  K\"ahler-Ricci flow  with  any initial K\"ahler metric in  $\mathcal{K}_{X}$   will converge to a  K\"ahler-Ricci soliton in
 $C^\infty$  in the sense of K\"ahler potentials.   Moreover, the convergence can be made
exponentially.

\end{theorem}

 We note that   that without loss of generality we may assume that
a K\"{a}hler-Ricci soliton $g_{KS}$ on $M$   is corresponding  to the above $X$ (cf. [TZhu1], [TZhu2]).
Theorem \ref{theorem-2} was first proved  by Tian and Zhu in [TZhu3] by  using an inequality of Moser-Trudinger type
established in [CTZ]\footnote{We need to add more details
about how to use the Moser-Trudinger typed inequality in general  case.}.    Here we will  modify arguments  in [TZhu5] in our general  case that $(M,J) $  admits a K\"ahler-Ricci soliton
 to give  an alternative proof of this theorem. This new
proof does not use such an  inequality of Moser-Trudinger type.   Moreover, in particular,  in case that $(M,J)$  admits  a K\"ahler-Einstein metric
 this new proof allows us to  avoid  to use a deep  result recentlly proved by Chen and Sun  in [CS]  for  the uniqueness of  K\"ahler-Einsteins in the sense of orbit space to give a self-contained proof to the main theorem in [TZhu5].

The organization of paper is as follows. In Section 1,
we discuss an upper bound of $\lambda(\cdot)$  in
general case-without any condition for $\mu(\cdot)$ and show that
the quantity $(2\pi)^{-n}[nV-N_X(c_1(M))]$ is an upper bound of
$\lambda(\cdot)$ in $\mathcal{K}_{X}$  (cf.  Proposition 1.4).   In Section 2, we
will summarize to give  some estimates for modified  Ricci potentials of
evolved K\"ahler metrics  along K\"ahler-Ricci flow
(cf. Proposition 2.3).
In Section 3, we prove Proposition 3.1 and  so  do Theorem 0.1.   Theorem 0.2 will be proved in Section 6.
  In Section 4,  we improve our key Lemma 3.2  in Section 3 independent of time $t$ (cf.  Proposition 4.2).   Section 5  is a discussion  about  an upper bound of
$\lambda(\cdot)$ in $\mathcal{K}_{Y}$ for a general  holomorphic vector field $Y\in
\eta_r(M)$. Section 7 is an appendix where we discuss the gradient estimate and Laplace estimate for the minimizers of Perelman's $W$-functional along the  K\"ahler-Ricci flow.

\vskip6mm
\section{ An upper bound of $\lambda(\cdot)$}

In this section, we  first review  Perelman's $W$-functional for triples
$(g,f,\tau)$ on a closed $m$-demensional Riemannian manifold $M$
(cf. [Pe], [TZhu5]).  Here  $g$ is a Riemannian metric, $f$ is a
smooth function and $\tau$ is a constant. In  our situation, we will
normalize volume of $g$ by
  \begin{equation}
\label{vol-1} \int_M dV_g\equiv\, V
 \end{equation}
and so we can fix  $\tau$ by $\frac{1}{2}$. Then the $W$-functional
depends only on a pair $(g,f)$ and it can be reexpressed as follows:
\begin{align}
\label{w-fun-1} W(g,f)=(2\pi)^{-m/2}\int_M[\frac{1}{2}(R(g)+ |\nabla
f|^2)+f] e^{-f}dV_g,\end{align}
 where $R(g)$ is a scalar curvature of $g$ and $(g,f)$  satisfies
 a normalization condition
\begin{equation}
\label{norm-1} \int_M e^{-f} dV_g\,= V.
\end{equation}
  Then  Perelman's entropy
$\lambda(g)$ is defined by
\begin{equation*}
\lambda(g)=\inf_{f}\{W(g,f)| ~(g,f)~\text{satisfies}~ (\ref{norm-1})
\}.
\end{equation*}

It is well known that $\lambda(g)$ can be attained by some smooth function $f$ (cf. [Ro]).
In fact, such a $f$ satisfies the Euler-Lagrange equation of $W(g,\cdot)$,
\begin{equation}\label{EulerLagrange}
\triangle f+f+\frac{1}{2}(R-|\nabla
f|^{2})=(2\pi)^{m/2}V^{-1}\lambda(g).
\end{equation}
Following Perelman's computation in [Pe], we can deduce the first
variation of $\lambda(g)$,
\begin{equation}\label{variation}
\delta\lambda(g)=-(2\pi)^{-m/2}\int_{M}<\delta g,
\mathrm{Ric}(g)-g+\nabla^{2}f>e^{-f}dV_{g},
\end{equation}
where $\mathrm{Ric}(g)$ denotes the Ricci tensor of $g$ and
$\nabla^{2} f$ is the Hessian of $f$. Hence, $g$ is a critical point
of $\lambda(\cdot)$ if and only if $g$ is a gradient shrinking
Ricci-soliton which satisfies
\begin{equation}\label{RicciSolitons}
\mathrm{Ric}(g)+\nabla^{2}f=g,
\end{equation}
where $f$ is  a  minimizer of $W(g,\cdot)$.
The following lemma  was proved in [TZhu5] for the uniqueness of
solutions  (\ref{EulerLagrange}) when $g$ is a gradient shrinking
Ricci soliton.

 \begin{lemma}\label{lemma-1}
If $g$ satisfies (\ref{RicciSolitons}) for some $f$, then any
solution of (\ref{EulerLagrange}) is equal to $f$  modulo a constant. Consequently,  a minimizer
of $W(g,\cdot)$ is unique if  the metric $g$ is a gradient shrinking
Ricci-soliton.
 Conversely, if $f$ is a function in (\ref{RicciSolitons}) for $g$, then $f$ satisfies (\ref{EulerLagrange}).
\end{lemma}

In case that $(M,J)$ is an $n$-dimensional  Fano manifold,   for any
K\"ahler metric $g$ in $2\pi c_1(M)$, (\ref{vol-1}) is equal to
  \begin{equation}
\label{vol-2} \int_M dV_g=\int_M \omega_g^n \,=\,(2\pi)^n\int_M
c_1(M)^n\,\equiv\, V.
 \end{equation}
Moreover,   (\ref{RicciSolitons}) becomes an equation for
K\"ahler-Ricci solitons,
$$
Ric(\omega_g)-\omega_{g}=\mathrm{L}_{X}\omega_{g},$$
  where $Ric(\omega_g)$ is a Ricci form of $g$ and  $L_X$ denotes the Lie derivative along a holomorphic vector
  field $X$ on $M$.  By
the uniqueness of K\"ahler-Ricci solitons [TZhu1], [TZhu2], we may
assume that   $X$   lies in    a
reductive Lie subalgebra $\eta_r(M)$ of $\eta(M)$ after a
holomorphism    transformation,     where  $\eta(M)$   consists of all
holomorphic vector
 fields on $M$.    Such a $X$ ( we call it an extremal  holomorphic vector
  field   for  K\"ahler-Ricci solitons ) can be determined as follows.

 Let  $\text{Aut}_r(M)$ be  a
 connected Lie subgroup of  automorphisms
group of $M$  generated by $\eta_r(M)$.   Let $K$ be a maximal compact subgroup of
$\text{Aut}_r(M)$. Without loss of generality, we may choose a $K$-invariant
 background  metric $g$   with its K\"ahler form $\omega_g$
in $2\pi c_1(M)$. In [TZhu2], as an obstruction to
K\"ahler-Ricci solitons, Tian and Zhu introduced a modified
Futaki-invariant $F_X(v)$ for any $X,v~\in~\eta(M)$ by
\begin{equation}\label{futaki-inv}
F_{X}(Z)=\int_{M}Z(h_g-\hat\theta_{X,\omega_g})e^{\hat\theta_{X,\omega_g}}\omega_g^{n},\quad\forall
~Z\in\eta(M),
\end{equation}
where $h_g$ is a Ricci potential of $g$ and $\hat\theta_{X,\omega_g}$  is  a real-valued potential of $X$  associated  to  $g$ defined by
$L_{X}\omega_g=\sqrt{-1}\partial\bar{\partial}\hat\theta_{X,\omega_g}$
with a normalization condition
\begin{equation}
\int_{M}\hat\theta_{X,\omega_g}   e^{h_g}\omega_g^{n}=0.
\end{equation}
It was showed
that there exists a unique   $X~\in ~\eta_r(M)$  such that
$$F_X(v)\equiv 0,~\forall~v\in~\eta_r(M).$$
Moreover, $F_X(v)\equiv 0$, for any $v\in~\eta(M)$ if $(M,J)$ admits  a
K\"ahler-Ricci soliton.

Let   $K_X$  be  an  one-parameter  compact subgroup  of holomorphisms transformation group  generated by   $X$.
We denote $\mathcal{K}_{X}$  to be a  class of $K_X$-invariant  K\"{a}hler metrics in $2\pi c_{1}(M)$.
Let  $\theta_{X,\omega_g}$  be  a real-valued potential of $X$  associated to  $g$
with a normalization condition
\begin{equation}
\int_{M}e^{\theta_{X,\omega_g}}\omega_g^{n}=\int_{M}\omega_g^{n}=V.
\end{equation}
Clearly,   $\theta_{X,\omega_g}=\hat\theta_{X,\omega_g}-c_X$  for some  constant $c_X$ which is independent of $g\in\mathcal{K}_{X} $ .

\begin{definition} For $g\in\mathcal{K}_{X}$, define $N_X(\omega_g)$
by
$$
N_X(\omega_g)=\int_{M}\theta_{X,\omega_g}e^{\theta_{X,\omega_g}}\omega_{g}^{n}.
$$
\end{definition}

By  Jensen's inequality,
it is easy to see
\begin{eqnarray*}
\begin{aligned}
&\frac{1}{V}\int_{M}(-\theta_{X,\omega_g})e^{\theta_{X,\omega_g}}\omega_{g}^{n}\\
&\le \text{log} \{\frac{1}{V}\int_{M}
e^{-\theta_{X,\omega_g}}e^{\theta_{X,\omega_g}}\omega_{g}^{n}\}=0.
\end{aligned}
\end{eqnarray*}
The equality holds if and only if $\theta_{X,\omega_g}=0$.
This shows that  $N_X(\omega_g)$ is  nonnegative and it is zero if and only if  the Futaki-invariant
vanishes [Fu].  Moreover, we have

\begin{lemma}\label{lemma-2}
 $N_X(\omega_g)$  is independent of choice of $g$ in $ \mathcal{K}_{X}$.
\end{lemma}

\begin{proof}

Choose a $K$-invariant K\"ahler form   $\omega$ in $2\pi c_1(M)$.
Then   for any K\"ahler metric $g$ in $ \mathcal{K}_{X}$ there
exists a K\"ahler potential $\varphi$ such that  the imaginary part of
 $X(\varphi)$ vanishes and  K\"ahler form of $g$ satisfies
$$\omega_g=\omega_\varphi=\omega+\sqrt{-1}\partial\bar{\partial}\varphi.$$
  Thus we suffice to prove
   $$N_X(\omega_\varphi)=N_X(\omega_{t\varphi}), ~~\forall ~t\in ~[0,1],$$
where
$\omega_{t\varphi}=\omega +t\frac{\sqrt{-1}}{2\pi}\partial\bar{\partial}\varphi.$
This  follows from
\begin{eqnarray*}
\begin{aligned}
\frac{dN_X(\omega_{t\varphi})}{dt}&=\int_{M}X(\varphi)e^{\theta_{X,\omega_{t\varphi}}}\omega_{t\varphi}^{n}+\int_{M}\theta_{X,\omega_{t\varphi}}
(X+\triangle)(\varphi)e^{\theta_{X,\omega_{t\varphi}}}\omega_{t\varphi}^{n}\\
&=\int_{M}X(\varphi)e^{\theta_{X,\omega_{t\varphi}}}\omega_{t\varphi}^{n}
-\int_{M}\nabla^{i} \varphi\nabla_{\bar{i}}\theta_{X,\omega_{t\varphi}}\omega_{t\varphi}^{n}\\
&=0.
\end{aligned}
\end{eqnarray*}
Here we have used  the fact
 $$\theta_{X,\omega_{t\varphi}}=\theta_{X,\omega_0}+tX(\varphi).$$
\end{proof}

By the above lemma, $N_X(\cdot)$ is an invariance  on
$\mathcal{K}_{X}$, which is independent of choice of $g$. For
simplicity,  we denote this invariance by $N_X(c_1(M))$. The
following  proposition  gives an upper bound of $\lambda(\cdot)$ in
$\mathcal{K}_{X}$ related to $N_X(c_1(M))$.

\begin{proposition}\label{upperbound}
$$\sup_{g\in\mathcal{K}_{X}}\lambda(g)\le(2\pi)^{-n}[nV-N_X(c_1(M))].$$
\end{proposition}

\begin{proof} Since $\lambda(g)\le W(g,-\theta_{X,\omega_{g}})$, we
suffice to prove
 \begin{equation}
 W(g,-\theta_{X,\omega_{g}})=(2\pi)^{-n}[nV-N_X(c_1(M))].
 \end{equation}
In fact, by using the facts $R(g)=2n+\Delta h_g$ and
$$\int_{M}(\Delta
\theta_{X,\omega_{g}}+|\nabla
\theta_{X,\omega_{g}}|^{2})e^{\theta_{X,\omega_{g}}}\omega_{g}^{n}=0,$$
we have
\begin{eqnarray*}
\begin{aligned}
&\int_{M}(R(g)+|\nabla
\theta_{X,\omega_{g}}|^{2})e^{\theta_{X,\omega_{g}}}\omega_{g}^{n}\\
&=2nV+\int_{M}(\Delta h_g-\Delta \theta_{X,\omega_{g}})e^{\theta_{X,\omega_{g}}}\omega_{g}^{n}\\
&=2nV-\int_{M}<\nabla(h_g-\theta_{X,\omega_{g}}),\nabla \theta_{X,\omega_{g}}>e^{\theta_{X,\omega_{g}}}\omega_{g}^{n}\\
&=2nV-2\int_{M}X(h_g-\theta_{X,\omega_{g}})e^{\theta_{X,\omega_{g}}}\omega_{g}^{n}\\
&=2nV-2e^{-c_X}F_{X}(X).
\end{aligned}
\end{eqnarray*}
In the last equality above, we  used  the relation  (\ref{futaki-inv}).  Since $X$ is extremal , we have
$$F_{X}(X)=0.$$
Thus by (\ref{w-fun-1}) for $f=-\theta_{X,\omega}$ together with Lemma
\ref{lemma-2}, one will get (1.11).
\end{proof}

In case that $M$ admits a  K\"{a}hler-Ricci soliton  $g_{KS}$, by
Lemma \ref{lemma-1},  a minimizer $f$ of $W(g_{KS},\cdot)$ in
$\mathcal{K}_{X}$ must be $-\theta_{X}$. Thus for any
$g\in\mathcal{K}_{X}$, by Proposition 1.4, we have
\begin{eqnarray*}
\begin{aligned}
\lambda(g_{KS})&=W(g_{KS},-\theta_{X})\\
               &=(2\pi)^{-n}[nV-N_X(c_1(M))]\ge\lambda(g).
               \end{aligned}
               \end{eqnarray*}
Therefore we get the following corollary.

\begin{corollary}\label{coro-1}
Suppose that $(M,J)$ admits a  K\"{a}hler-Ricci soliton $g_{KS}$. Then
$g_{KS}$ is a global maximizer of $\lambda(\cdot)$ in
$\mathcal{K}_{X}$ and
\begin{equation}
\lambda(g_{KS})=(2\pi)^{-n}[nV-N_X(c_1(M))].
\end{equation}

\end{corollary}

\begin{remark}  Corollary \ref{coro-1} implies that a K\"ahler-Einstein metric  is a global maximizer of $\lambda(\cdot)$  in $2\pi c_1(M)$ even with varying complex structures
and  supremum of $\lambda(\cdot)$
is  $(2\pi)^{-n}nV$  since  $N_X(c_1(M))=0$.    Note that $N_X(c_1(M))>0$ if the
Futaki-invariant does not vanish. Thus Corollary \ref{coro-1} also
implies that the supremum of $\lambda(\cdot)$  in case that $(M,J)$
admits a K\"{a}hler-Ricci soliton is strictly less than one in case
that $(M,J)$ admits a K\"{a}hler-Einstein metric.
\end{remark}

\vskip6mm

\section{Estimates for modified Ricci potentials }

In this section, we summarize  some apriori estimates for modified
Ricci potentials of evolved K\"ahler metrics along  K\"ahler-Ricci
flow.  Some similar estimates have been also discussed in [TZhu3] and
[PSSW],  we refer the readers to those  two papers.  We consider the following
(normalized) K\"ahler-Ricci flow:
\begin{equation}\label{KRF}
\frac{\partial g(t,\cdot)}{\partial t}\,=\,-{\rm Ric}(g(t,\cdot)) +
g(t,\cdot),~~~g(0)=g,
\end{equation}
where  K\"ahler form of $g$ is in $2\pi c_1(M)$. It was proved in
[Ca] that (\ref{KRF}) has a global solution $g_t=g(t,\cdot)$ for all
time $t>0$. For simplicity, we denote by $(g_t;g)$ a solution of
(\ref{KRF})  with initial metric $g$. Since the flow preserves the
K\"ahler class, we may write K\"ahler form of $g_t$ as
$$\omega_{\phi}\,=\,\omega_g+\sqrt{-1}\partial\overline\partial \phi$$
for some K\"aher potential $\phi=\phi_t$.

Let $X\in\eta_r(M)$ be  the extremal  holomorphic vector field on
$M$ as in Section 1 and  $\sigma_t=\exp\{tX\}$  an  one-parameter
subgroup  generated by $X$.  Let $\phi'=\phi_{\sigma_t}$ be
corresponding K\"ahler potentials of
$\sigma_t^\star\omega_{\phi_t}$.  Then $\omega_{\phi'}$ will
satisfy a modified K\"ahler-Ricci flow,
 \begin{equation}\label{equation-2}
 \frac
{\partial}{\partial t}\omega_{\phi'}=-\text
{Ric}(\omega_{\phi'})+\omega_{\phi'}+L_X\omega_{\phi'}.
 \end{equation}
Equation (\ref{equation-2}) is equivalent to the following
Monge-Amp\'ere flow for $\phi'$ (modulo a constant),
 \begin{equation}\label{equation-3}
\frac {\partial\phi'}{\partial t}=\log \frac
 {\omega_{\phi'}^n}{\omega_g^n}+\phi'+\theta_{X,\omega_{\phi'}}
  -h_g,~~~~ \phi'(0,\cdot)=c,
\end{equation}
where $c$ is a constant and all  K\"ahler potentials $\phi'=\phi_t'=\phi'(t,\cdot)$  are in
a space given by
$$\mathcal{P}_{X}(M,\omega)=\{\varphi\in
C^{\infty}(M)|~\omega_{\varphi}=\omega+\sqrt{-1}\partial\overline{\partial}\varphi>0,\quad
\mathrm{Im}(X(\varphi))=0\}.$$
 
 By using the maximum principle to (\ref{equation-2}) or
 (\ref{equation-3}),  we get
 \begin{equation}\label{relation-1}
  h_{\phi'}-\theta_{X,\omega_{\phi'}}=-\frac {\partial}{\partial
t}\phi'+c_t,
\end{equation}
 for some constants $c_t$.  Here $h_{\phi'}$ are Ricci potentials of
 $\omega_{\phi'}$ which are normalized by
\begin{equation}\label{normalization-3}
  \int_M e^{h_{\phi'}}\omega_{\phi'}^n=V.
  \end{equation}
The following estimates are due to G. Perelman. We refer the  readers to
[ST] for their proof.

\begin{lemma}\label{lem:perelman-1}  There are constants $c$ and $C$ depending only on the initial metric $g$ such that (a)
$\text{diam}(M,\omega_{\phi'})\le C$; (b)
$\text{vol}(B_r(p),\omega_{\phi'} )\ge c r^{2n}$; (c) $
\|h_{\phi'}\|_{C^0(M)}\le C$; (d) $\|\nabla
h_{\phi'}\|_{\omega_{\phi'}} \le C$; (e) $\|\Delta h_{\phi'}
\|_{C^0(M)} \le C$.
 \end{lemma}

  Recall that the modified Mabuchi's K-energy  $\mu(\cdot)$ is
  defined in $\mathcal{P}_{X}(M,\omega)$ by
\begin{eqnarray*}
\begin{aligned}
\mu (\varphi)
 =~&-\frac {n}{ V}\int_0^1\int_M
\dot\psi[\text
{Ric}(\omega_{\psi})-\omega_{\psi}
-\sqrt{-1}\partial\overline\partial
\theta_{X,\omega_\psi}\\
 &  +\sqrt{-}\overline\partial
(h_{\omega_{\psi}}-\theta_{X,\omega_\psi})\wedge\partial\theta_{X,\omega_\psi})
]\wedge e^{\theta_{X,\omega_\psi}}\omega_{\psi}^{n-1} \wedge dt,
\end{aligned}
\end{eqnarray*}
 where $\psi=\psi_t$ $(0\le t\le 1)$ is a path connecting $0$ to $\varphi$ in $\mathcal{P}_{X}(M,\omega)$.
 If $X=0$, then $\mu_{\omega_g}(\phi)$ is nothing but
Mabuchi's $K$-energy [Ma].
 Then by (\ref{equation-2}), we have
\begin{equation}\label{monoto-2}
\frac {d\mu(\phi')}{dt} =-\frac {1}{V}\int_M\|\overline\partial\frac
{\partial\phi'}{\partial t}\|_{\omega_{\phi'}}^2
 e^{\theta_{X,\omega_{\phi'}}}(\omega_{\phi'})^n\le 0.
\end{equation}
This implies that $\mu(\phi')$ is uniformly bounded if $\mu(\cdot)$
is bounded from below in $\mathcal{P}_{X}(M,\omega)$.

 Let
$u_{X,\phi'}=u_{X,\omega_{g_t'}}=h_{\phi'}-\theta_{X,\omega_{\phi'}}$.
Then

\begin{lemma}\label{lemma} There exists a uniform $C$ such that
$$\|\nabla u_{X,\phi'}\|_{\omega_{\phi'}}\le C.$$
\end{lemma}

\begin{proof} First we note that $\theta_{X,\omega_{\phi'}}$ is uniformly bounded in
$\mathcal{P}_{X}(M,\omega)$ (cf. [Zhu1], [ZZ]). Then by (c) of Lemma
\ref{lem:perelman-1}, we have
$$\|u_{X,\phi'}\|_{C^0}=\|u_{X,\omega_{g_{s}'}}\|_{C^0}\le C, ~\forall ~s>0$$
 for some uniform constant $C$. Now we consider the flow
 (\ref{equation-3}) with zero as an initial K\"ahler potential and
 the background K\"ahler form $\omega_g$  replaced by $\omega_{g_s'}$. By an
 estimate in Lemma 4.3 in [CTZ], we see
 $$t\|\nabla u_{X,\omega_{g_{s+t}'}}\|_{\omega_{g_{s+t}'}}^2\le
 e^{2t}\|u_{X,\omega_{g_{s}'}}\|_{C^0},~~\forall~t>0.$$
In particular, we get
 $$ \|\nabla u_{X,\omega_{g_{s+t}'}}\|_{\omega_{g_{s+t}'}}^2\le C',~~\forall~t\in
 [1,2].$$
Since the above estimate is independent of $s$, we conclude that the
lemma is true.

\end{proof}

Now we begin to prove the main result in this section.

\begin{proposition}\label{Limit}
Suppose that  $\mu(\cdot)$  is bounded from below in $\mathcal{P}_{X}(M,\omega)$. Then we have:\\
(a) $\lim_{t\rightarrow\infty}||u_{X,\phi'}||_{C^{0}}=0$;\\
(b) $\lim_{t\rightarrow\infty}||\nabla u_{X,\phi'}||_{\omega_{\phi'}}=0$;\\
(c) $\lim_{t\rightarrow\infty}||\triangle u_{X,\phi'}||_{C^{0}}=0$.
\end{proposition}

\begin{proof} Let $H(t)=\int_{M}|\nabla u_{X,\omega_{g_t}'}|^{2}e^{\theta_{X,\omega_{g_t}'}}\omega_{g_t'}^{n}$.  Then by
(\ref{monoto-2}),  one sees that there exists a sequence of
$t_{i}\in[i,i+1]$ such that
 \begin{equation*}
 \lim_{i\rightarrow\infty}H(t_{i})=0.
 \end{equation*}
 Thus by using a differential inequality
 \begin{equation*}
 \frac{dH(t)}{dt}\leq CH(t),
 \end{equation*}
 where $C$ is a uniform constant (cf. [PSSW]), we get
 \begin{equation}\label{nabla-esti}
 \lim_{t\rightarrow\infty}\int_{M}|\nabla u_{X,\omega_{g_t'}}|_{g_t'}^{2}e^{\theta_{X,\omega_{g_t'}}}\omega_{g_t'}^{n}=0.
 \end{equation}

 Let
 \begin{equation*}
 \tilde{u}_{t}=u_{X,\omega_{g_t'}}
 -\frac{1}{V}\int_Mu_{X,\omega_{g_t'}}e^{h_t'}\omega_{g_t'}^{n},
 \end{equation*}
 where $h_t'=h_{\phi'(t,\cdot)}$. Then by using the weighted Poincar\'{e} inequality in [TZhu3]
 together with (c) of Lemma \ref{lem:perelman-1}, we
 obtain
 from (\ref{nabla-esti}),
 \begin{equation*}
 \int_{M}\tilde{u}_{t}^{2}e^{h_t'}\omega_{g_t'}^{n}\leq\int_{M}|\nabla u_{X,\omega_{g_t'}}|_{g_t'}^{2}e^{h_t'}\omega_{g_t'}^{n}\to 0,~~\text{as }~t\to\infty.
 \end{equation*}
 Consequently, we derive
 \begin{equation}\label{L2-esti}
\lim_{t\rightarrow\infty}\int_{M}\tilde{u}_{t}^{2}\omega_{g_t'}^{n}=0.
 \end{equation}

 We claim
 \begin{equation*}
 \lim_{t\rightarrow\infty}\|\tilde{u}_{t}\|_{C^{0}}=0.
 \end{equation*}
The claim immediately implies (a) of Proposition \ref{Limit} by the
normalization conditions
\begin{equation*}
\int_{M}e^{\theta_{X,\omega_{g_t'}}}\omega_{g_t'}^{n}=\int_{M}e^{h_t'}\omega_{g_t'}^{n}=V.
\end{equation*}
To prove the claim, we need to use an inequality
\begin{equation}\label{harnack}
\|\tilde{u}_{t}\|_{C^{0}}^{n+1}\leq C\|\nabla
u_{X,\omega_{g_t'}}\|^{n}_{g_t'}[\int_{M}\tilde{u}_{t}^{2}\omega_{g_t'}^{n}]^{\frac{1}{2}}.
\end{equation}
(\ref{harnack}) can be proved by using the non-collapsing estimate
(b) in Lemma 2.1 (cf. [PSSW], [Zhu2]). Thus by Lemma 2.2 and
(\ref{L2-esti}), the claim is proved.

By (a)  we can  show that after a suitable choice of constant $c$
in the flow (\ref{equation-3}) it holds
 $$\lim_{t\to \infty}\|\frac {\partial}{\partial
t}\phi'\|_{C^0}=0.$$
   In fact under the assumption of lower bound
of modified $K$-energy, one can choose such a  $c$ (cf.
[TZhu3])
 such that
$$\lim_{t\to\infty}\int_M\frac {\partial}{\partial
t}\phi' e^{\theta_{X,\omega_{\phi'}}}\omega_{\phi'}^{n}=0.$$
 Then by (\ref{relation-1}),
  we will get the conclusion. On the other hand, by
 Lemma \ref{lemma} and (d) of Lemma 2.1, we have
$$\sup_{t\in[0,\infty)}\|X\|_{g_t'}<C$$
 for some uniform constant $C$. Therefore, by using the following
 lemma we prove  (b) and (c).

\end{proof}

\begin{lemma}{([PSSW])}
There exist $\delta, K>0$ depending only on $n$ and the constant
$C_{X}=\sup_{t\in[0,\infty)}\|X\|_{g_t'}$ with the following
property. For $\epsilon$ with $0<\epsilon\leq \delta$ and any
$t_{0}>0$, if
\begin{equation*}
\|\frac{\partial\phi'}{\partial t}\|_{C^{0}}(t_0)\leq \epsilon,
\end{equation*}
then
\begin{equation*}
\|\nabla u_{X,\omega_{g_{t_0+2}'}}\|^2_{g_{t_0+2}'}+\|\Delta
u_{X,\omega_{g_{t_0+2}'}}\|_{C^{0}}\leq K\epsilon.
\end{equation*}
\end{lemma}

\vskip6mm
\section{Proof of Theorem \ref{theorem-1}}

According to [TZhu5],  an energy level $L(g)$ of entropy
$\lambda(\cdot)$ along  K\"{a}hler-Ricci flow $(g_t;g)$  is
defined by
\begin{equation*}
L(g)=\lim_{t\rightarrow\infty}\lambda(g_t).
\end{equation*}

By the monotonicity of $\lambda(g_t)$, we see that $L(g)$ exists and
it is finite. In this section, our goal is to prove

\begin{proposition}\label{prop-2}
Suppose that the modified Mabuchi's K-energy is bounded from below
in $\mathcal{K}_{X}$. Then for any $g\in\mathcal{K}_{X}$.
\begin{equation}\label{level}
L(g)=(2\pi)^{-n}(nV-N_X(c_1(M)).
\end{equation}

\end{proposition}

The above proposition  shows that the energy level $L(g)$ of
entropy $\lambda(\cdot)$ does not depend on the initial K\"ahler
metric $g\in\mathcal{K}_{X}$. Thus by using the K\"ahler-Ricci
flow $(g_t;g)$ for any K\"ahler metric $g\in\mathcal{K}_{X}$ and
the monotonicity of $\lambda(g_t)$, we will get Theorem
\ref{theorem-1}.

To prove Proposition \ref{prop-2}, we need the following  key lemma.

\begin{lemma}\label{keylemma}  Let $f_t$ be  a   minimizer of $W(g_t,\cdot)$-functional associated  evolved K\"ahler metric $g_t$  of  (2.1) at time $t$ 
and $h_t$ a   Ricci potential of $g_t$ which satisfying the normalization (\ref{normalization-3}) . Then there exists a sequence of $t_i\in [i,i+1]$ such that

\noindent (a) $\lim_{t_i\to \infty}\|\Delta (f_{t_i}+h_{t_i})
\|_{L^2(M,\omega_{g_{t_i}})}\,=\,0$;

\noindent (b) $\lim_{t_i\to \infty}\|\nabla (f_{t_i}+h_{t_i})
\|_{L^2(M,\omega_{g_{t_i}})}\,=\,0$;

\noindent (c) $\lim_{t_i\to \infty} \|f_{t_i}+h_{t_i}\|_{C^0}\,=\,
0$.
\end{lemma}

\begin{proof}  Lemma \ref{keylemma} is a generalization of  Proposition 4.4  in
[TZhu5]. We will follow the argument there.    First  by (\ref{variation}),   it is easy to see that
\begin{eqnarray*}
\begin{aligned}
\frac{d}{dt}\lambda(g_t)=(2\pi)^{-n}\int_{M}|\mathrm{Ric}(g_t)-g_t+
\nabla^{2}f_t|_{g_t}^{2}e^{-f_t}\omega_{g_t}^{n}.
\end{aligned}
\end{eqnarray*}
It follows  that
\begin{equation*}
\frac{d}{dt}\lambda(g_t))\geq
(2\pi)^{-n}\frac{1}{2n}\int_{M}|\triangle(h_t+f_t)|^{2}e^{-f(t)}\omega_{g_t}^{n}.
\end{equation*}
Since $\lambda(g_t)\leq W(g_t,0)=(2\pi)^{-n}nV$  are uniformly
bounded, we see that there exists a sequence of $t_{i}\in [i,i+1]$
such that
\begin{equation*}
\lim_{i\rightarrow\infty}\int_{M}|\triangle(h_{t_{i}}+f_{t_{i}})|^{2}e^{-f_{t_{i}}}\omega_{g_{t_{i}}}^{n}=0.
\end{equation*}
Note that $f_t$ is uniformly bounded [TZhu5]. Hence we see that that
(a) of the lemma is true. By (a), we also get
\begin{eqnarray}
 \begin{aligned}
 &\lim_{t_i\to \infty}\|\nabla
(f_{t_i}+h_{t_i})
\|_{L^2(M,\omega_{g_{t_i}})}\\
&\le \lim_{t_i\to \infty}
\int_M|f_{t_i}+h_{t_i}||\triangle(f_{t_i}+h_{t_i})|\omega_{g_{t_{i}}}^{n}\\
&\le C\lim_{t_i\to \infty}\|\Delta (f_{t_i}+h_{t_i})
\|_{L^2(M,\omega_{g_{t_i}})}=0.
\end{aligned}
\end{eqnarray}
This proves (b) of the lemma. It remains to prove (c).

 Let $q_t=f_t+h_t$.  Then 
  \begin{eqnarray}\label{equation-4}
 \begin{aligned}
-\Delta {q_t}&=-\Delta f_t-\Delta
h_t\\
&=f_t+\frac{1}{2}(R-|\nabla
f_t|^{2})-(2\pi)^{2n}V^{-1}\lambda(g_t)-\Delta h_t\leq C.
\end{aligned}
 \end{eqnarray}
Define
\begin{equation*}
\tilde{q_t}=q_t-\frac{1}{V}\int_{M}q_te^{h_t}\omega_{g_t}^{n}.
\end{equation*}
By using the weighted Poincar\'{e} inequality (cf. [TZhu3]), we have
\begin{equation*}
\int_{M}\tilde{q_t}^{2}e^{h_t}\omega_{g_t}^{n}\leq \int_{M}|\nabla
q_t|^{2}e^{h_t}\omega_{g_t}^{n}.
\end{equation*}
It follows by (b),
 \begin{equation}
\lim_{i\rightarrow\infty}\int_{M} \tilde
q_{t_i}^{2}\omega_{g_{t_i}}^{n}=0.
\end{equation}
Hence,  following an argument in the proof of  Proposition 4.4  in
[TZhu5], we will get estimates
\begin{equation} \label{estimate-3}
\|\tilde{q_{t_i}}^{+}\|_{C^0}\leq C \|\tilde
q_{t_i}\|_{L^2(M,\omega_{g_{t_i}})}\to 0,~~\text{as}~i\to \infty,
\end{equation}
 and 
\begin{equation}\label{estimate-4}
\lim_{i\rightarrow\infty}\int_{M}\tilde{q_{t_{i}}}^{-}
\omega_{g_{t_{i}}}^{n}=0,
\end{equation}
where $q_t^{+}=\max\{q_t,0\}$ and $q_t^{-}=\min\{q_t,0\}$. 
Consequently,  we derive
$$\int_{M}\tilde{q_{t_i}}e^{-f_{t_{i}}}\omega_{g_{t_{i}}}^{n}=0.$$
This implies 
\begin{equation}\label{estimate-4}
\lim_{i\rightarrow\infty}\int_{M}q_{t_i}e^{-f_{t_i}}\omega^{n}_{g_{t_{i}}}=0
\end{equation}
 according to  the normalization  $\int_{M}e^{-f_t}\omega^{n}_{g_t}=\int_{M}e^{h_t}\omega^{n}_{g_t}=V$.

Next we improve that 
\begin{equation}\label{estimate-5}
\lim_{i\rightarrow\infty}|q_{t_i}|=0.
\end{equation}
 Let $u_t=e^{-\frac{f_t}{2}}-e^{\frac{h_t}{2}}$. We claim
\begin{equation}\label{claim}
\lim_{i\rightarrow\infty}\|u_{t_i}\|_{L^{2}(M,\omega_{g_{t_i}})}=0.
\end{equation}
In fact,  by Jensen's inequality and (\ref{estimate-4}), one sees
\begin{eqnarray*}
\begin{aligned}
\frac{1}{V}\int_{M}e^{-\frac{f_{t_i}}{2}}e^{\frac{h_{t_i}}{2}}\omega_{g_{t_i}}^{n}
&=\frac{1}{V}\int_{M}e^{\frac{f_{t_i}+h_{t_i}}{2}}e^{-f_{t_i}}\omega_{g_{t_i}}^{n}\\
&\geq e^{\frac{1}{2V}\int_{M}(f_{t_i}+h_{t_i)}e^{-f_{t_i}}\omega_{g_{t_i}}^{n}}\to
1, ~\text{as}~i\to \infty.
\end{aligned}
\end{eqnarray*}
 On the other hand,
\begin{equation*}
\int_{M}e^{-\frac{f_t}{2}}e^{\frac{h_t}{2}}\omega_{g_t}^{n}\leq(\int_{M}e^{-f_t}\omega_{g_t}^{n})^{\frac{1}{2}}(\int_{M}e^{h_t}\omega_{g_t}^{n})^{\frac{1}{2}}=V.
\end{equation*}
Hence
\begin{equation*}
\lim_{i\rightarrow\infty}\int_{M}e^{-\frac{f_{t_{i}}}{2}}e^{\frac{h_{t_{i}}}{2}}\omega_{g_{t_{i}}}^{n}=V.
\end{equation*}
It follows
\begin{equation*}
\lim_{i\rightarrow\infty}\int_{M}u_{t_i}^{2}\omega_{g_{t_{i}}}^{n}
=2V-2\lim_{i\rightarrow\infty}\int_{M}e^{-\frac{f_{t_{i}}}{2}}e^{\frac{h_{t_{i}}}{2}}\omega_{g_{t_{i}}}^{n}=0.
\end{equation*}
This completes the proof of claim.

Since equation (\ref{EulerLagrange}) is equivalent to
\begin{equation}\label{EulerLagrange-2}
\Delta v_t-\frac{1}{2}f_tv_t-\frac{1}{4}R(g_t)v_t=\frac{1}{2V}(2\pi)^{n}\lambda(g_t)v_t,
\end{equation}
 where $v_t=e^{\frac{-f_t}{2}}$, by Lemma
\ref{lem:perelman-1}, it is easy to see
\begin{equation*}
|\Delta u_t|\leq C.
\end{equation*}
 Then by the standard Moser's iteration, we
get from (\ref{claim}),
\begin{equation*}
\|u_{t_i}\|_{C^0}\le C\|u_{t_i}\|_{L^2(M,\omega_{g_{t_i}})}\to
0,~\text{as}~i\to\infty.
\end{equation*}
This implies (\ref{estimate-5}),  so we obtain (c) of the lemma.

\end{proof}

\begin{proof}[Proof of Proposition \ref{prop-2}]
Note that $\frac{R(g_t)}{2}=n+\frac{1}{2}\Delta h_t$, where
$\Delta$ is the Beltrima-Laplacian operator associated to
the Riemannian metric $g_t$.
Then
\begin{equation*}
\int_{M}\frac{1}{2}(R(g_t)+|\nabla
f_t|^{2})e^{-f_t}dV_{g_t}=nV+\frac{1}{2}\int_{M}\Delta
(f_t+h_t)e^{-f_t}dV_{g_t}.
\end{equation*}
Thus by (a) of Lemma \ref{keylemma}, one sees that there exists a
sequence of time $t_i$ such that
\begin{equation}\label{estimate-6}
\lim_{i\rightarrow\infty}\int_{M}\frac{1}{2}(R(g_{t_i})+|\nabla
f_{t_i}|^{2})e^{-f_{t_i}}dV_{g_{t_{i}}}=nV.
\end{equation}
On the other hand,  since the modified Mabuchi's K-energy is
bounded from below,  we see that  (a) of Proposition 2.3 is true.
Then by (c) of Lemma \ref{keylemma},  it follows
 \begin{equation}
 \lim_{i\to\infty}\|f_{t_i}+\theta_{X,\omega_{g_{t_i}}}\|_{C^0}=0.
\end{equation}
Here we used a fact $\sigma_t^\star\theta_{X,\omega_{g_{t}}}=\theta_{X,\omega_{g_{t}'}}$ since $X$ lies in the center of $\eta_r(M)$ [TZhu1].
Hence 
\begin{eqnarray}\label{estimate-7}
\begin{aligned}
&\lim_{i\to\infty}\int_{M}f_{t_i}e^{-f_{t_i}}dV_{g_{t_i}}\\
&=-\lim_{i\to\infty}
\int_{M}\theta_{X,\omega_{g_{t_i}}}e^{\theta_{X,\omega_{g_{t_i}}}}\omega^{n}_{g_{t_i}}=-N_X(c_1(M)).
\end{aligned}
\end{eqnarray}
By combining (\ref{estimate-6}) and (\ref{estimate-7}),   we get
\begin{eqnarray*}
\begin{aligned}
&\lim_{i\to\infty}
\lambda(g_{t_i})=\lim_{i\to\infty}\int_{M}[\frac{1}{2}(R(g_{t_i})+|\nabla
f_{t_i}|^{2})+f_{t_i}]e^{-f_{t_i}}dV_{g_{t_{i}}}\\
&=nV-N_X(c_1(M))
\end{aligned}
\end{eqnarray*}
Therefore, by using the monotonicity of $\lambda(g_t)$ along the
flow $(g_t;g)$, we  obtain (\ref{level}).

\end{proof}

It was showed in [TZhu4] that a K\"ahler-Ricci soliton is a local maximizer of $\lambda(\cdot)$  in  the K\"ahler  class $2\pi c_1(M)$.   Together with  Corollary  \ref{coro-1},
 one may guess that  a K\"ahler-Ricci soliton is a  global maximizer of $\lambda(\cdot)$.
  More general,   according to Theorem \ref{theorem-1} ,  we propose the following conjecture.

\begin{conjecture}\label{conjecture} Suppose that the modified Mabuchi's K-energy is bounded from below.
Then
\begin{equation*}
\sup_{\omega_{g'}\in 2\pi c_1(M)} \lambda({g'})=(2\pi)^{-n}[nV-N_X(c_1(M))].
\end{equation*}
\end{conjecture}

\vskip6mm
\section{Improvement of  Lemma 3.2}

In this section, we use Perelman's backward heat flow to improve estimate (c) in Lemma \ref{keylemma} independent of $t$.  Moreover, we show the gradient estimate of $f_t+h_t$ also holds.
Although Lemma  \ref{keylemma}  is sufficient to be applied  to  prove   Theorem 0.1 and Theorem 0.2,  results of this section are independent of interests.  We hope that these results will have
 applications  in the future.

Fix any $t_0\geq 1$.  We consider   a backward heat equation in $t\in[t_0-1,t_0]$,
\begin{equation}\label{backflow}
\frac{\partial}{\partial t}f_{t_0}(t)=-\triangle f_{t_0}(t)+|\nabla f_{t_0}(t)|^2-\triangle h_t
\end{equation}
with an initial $f_{t_0}(t_0)=f_{t_0}$.   Clearly,  the equation  preserves the normalizing condition $\frac{1}{V}\int_Me^{-f_{t_0}(t)}\omega_{g_t}^n=1$.
  Moreover, by the maximum  principle, we have
\begin{equation}\label{equation-4.1}
\|f_{t_0}(t)\|_{C^0}\leq C(g),\hspace{0.3cm}\forall~ t\in[t_0-1,t_0],
\end{equation}
since $\triangle h_t$ are uniformly bounded.   Here the constant $C(g)$ depends only on  the initial metric $g$  of (2.1).

Similarly  to (\ref{variation}),  we can compute
\begin{eqnarray}\label{variation2}
\begin{aligned}
&\frac{d}{dt}W(g_t,f_{t_0})\\
&=(2\pi)^{-n} \int_M(\|\partial\overline\partial (h_t+f_{t_0}(t))\|^{2}+\|\partial\partial f_{t_0}(t)\|^2 )e^{-f_{t_0}(t)}\omega_{_{g_t}}^n.
\end{aligned}
\end{eqnarray}
By   using (\ref{variation2}),   we  want to prove

\begin{lemma}\label{lemma-4.1}
\begin{equation}\label{estimate-4.3}
\|f_t+h_t-c_t\|_{L^2(M, g_t)}\rightarrow 0, ~\text{as} ~t\rightarrow\infty,
\end{equation}
where  $c_t=\frac{1}{V}\int_M(f_t+h_t)e^{h_t}\omega_{g_t}^n$.

\end{lemma}

\begin{proof}  First by   (\ref{variation2}),  one sees
\begin{eqnarray}
\lambda(g_{t_0})-\lambda(g_{t_0-1})&\geq& W(g_{t_0},f_{t_0}(t_0))   -W(g_{t_0-1},f_{t_0}(t_0-1))\nonumber\\
&\geq&(2\pi)^{-n}\frac{1}{2n}\int_{t_0-1}^{t_0}\int_M|\triangle(f_{t_0}(t)+h_t)|^2e^{-f_{t_0}(t)}\omega_{g_t}^ndt.\nonumber
\end{eqnarray}
It follows
\begin{equation}\nonumber
\int_{t_0-1}^{t_0}\int_M|\triangle(f_{t_0}(t)+h_t)|^2\omega^n_{g_t}dt\rightarrow 0, ~\text{as } ~t_0\rightarrow\infty.
\end{equation}
Thus by  using the weighted Poincar\'e inequality  as in (3.4) in last section ,  we  will get
\begin{eqnarray}\label{estimate-4.4}
\begin{aligned}
&\int_{t_0-1}^{t_0} dt \int_M(f_{t_0}(t)+h_t-c_{t_0}(t))^2\omega^n_{g_t}\\
&\leq C(g_0) [\int_{t_0-1}^{t_0} dt \int_M|\triangle(f_{t_0}(t)+h_t)|^2\omega^n_{g_t}]^{1/2}\to 0,  ~\text{as}  ~t_0 \to\infty,
\end{aligned}
\end{eqnarray}
 where $c_{t_0}(t)=\frac{1}{V}\int_M(f_{t_0}(t)+h_t)e^{h_t}\omega_{g_t}^n.$

Next,  since  $\frac{dh_t}{dt}=\Delta h_t+h_t-a_t$,  where $a_t=\frac{1}{V}\int_Mh_te^{h_t}\omega^n_{g_t}$,  by a straightforward calculation,  we see
\begin{equation}
\begin{aligned}
&\frac{d}{dt}\int_M(h_t+f_{t_0}-c_{t_0}(t))^2\omega^n_{g_t}\\
&=\int_M[2(h_t+f_{t_0}(t)-c_{t_0}(t))(\triangle f_{t_0}-|\nabla f_{t_0}(t)|^2+h_t-a_t-\frac{dc_{t_0}}{dt})
\nonumber\\
&-(h_t+f_{t_0}(t)-c_{t_0}(t))^2\triangle h_t] \omega^n_{g_t}\nonumber\\
\end{aligned}
\end{equation}
Then by  Lemma 3.2,  we get
\begin{equation}
\begin{aligned}
&|\frac{d}{dt}\int_M(h_t+f_{t_0}(t)-c_{t_0})^2 \omega^n_{g_t}|\\
&\leq C+C\int_M\big(|\triangle f_{t_0}(t)|+|\nabla f_{t_0}(t)|^2+|\frac{dc_{t_0}(t)}{dt}|\big) \omega^n_{g_t}\nonumber\\
&\leq C+C\int_M\big(|\nabla\bar{\nabla}(f_{t_0}(t)+h_t)|^2+|\frac{dc_{t_0}(t)}{dt}|\big)  \omega^n_{g_t}.\nonumber
\end{aligned}
\end{equation}
Notice that
\begin{eqnarray*}
\frac{d c_{t_0}}{dt}=\frac{1}{V}\int_M[\triangle f_{t_0}(t)-|\nabla f_{t_0}(t)|^2-(h_t+f_{t_0}(t))(h_t+a_t)]e^{h_t} \omega^n_{g_t}.
\end{eqnarray*}
We can also estimate
$$|\frac{dc_{t_0}}{dt}|\leq C+C\int_M|\nabla\bar{\nabla}(f_{t_0}(t)+h_t)|^2 \omega^n_{g_t}.$$
Hence we derive
\begin{equation}\label{estimate-4.5}
\big|\frac{d}{dt}\int_M(f_{t_0}(t)+h_t-c_{t_0}(t))^2dv\big|\leq C+C\int_M|\nabla\bar{\nabla}(f_{t_0}(t)+h_t)|^2 \omega^n_{g_t}.
\end{equation}
Therefore,  according to 
  \begin{equation*}
\begin{aligned}
&\int_{t_0-1}^{t_0} dt \int_M|\nabla\bar{\nabla}(f_{t_0}(t)+h_t)|^2 e^{-f_{t_0}(t)} \omega^n_{g_t}\\
&\le  (2\pi)^n (\lambda(g_{t_0})-\lambda(g_{t_0-1}))
\to 0, ~\text{as} ~t_0\rightarrow\infty,
\end{aligned}
\end{equation*}
 (\ref{estimate-4.4}) and  (\ref{estimate-4.5}) will implies
$$\|f_{t_0}(t)+h_t-c_{t_0}(t)\|_{L^2(g_t,M)}\rightarrow 0, ~\text{as} ~t_0\rightarrow\infty, ~\forall ~t\in [t_0-1,t_0].$$
 Consequently,   we get
(\ref{estimate-4.3}).
\end{proof}

\begin{proposition}
\begin{equation}\label{proposition-zhang}
\|f_t+h_t\|_{C^0}+\|\nabla(f_t+h_t)\|_{g_{t}}\rightarrow 0,
~\text{ as} ~t\rightarrow\infty. \end{equation}

\end{proposition}

\begin{proof} With the help of  Lemma \ref{lemma-4.1},  by using  same  argument in  the proof of (c) in Lemma 3.2, we
can prove that
 \begin{equation}\label{estimate-C^0}
\|f_t+h_t\|_{C^0}\rightarrow 0,
~\text{ as} ~t\rightarrow\infty. \end{equation}
So we suffice to prove
\begin{equation}\label{gradients-4.1}
\|\nabla(f_t+h_t)\|_{g_{t}}\rightarrow 0,
~\text{ as} ~t\rightarrow\infty. \end{equation}
We will use the Moser's iteration to obtain  (\ref{gradients-4.1} ) as in lemma 7.2 in Appendix.   We note by (\ref{estimate-C^0}) and  Theorem 7.1   that
\begin{equation}\label{estimate-4.6}
\int_{M}|\nabla q_t|^{2}\omega_{g_t}^n=|\int_{M} -\Delta (f_t+h_t) (f_t+h_t)\omega_{g_t}^n|\leq   C \|f_t+h_t\|_{C^{0}}\to 0,\end{equation}
where $q_t=  f_t+h_t$  satisfies an equation
$$\triangle q_t=\frac{1}{2}(|\nabla f_t|^2-2f_t+\triangle h_t)+(2\pi)^{n}V^{-1}\lambda(g_t)-n.$$

Let $w_t=|\nabla q_t|^2$.   Then by  the Bochner formula, we have
\begin{eqnarray}
\triangle w_t&=&|\nabla\nabla q_t|^2+|\nabla\bar{\nabla}q_t|^2+\nabla_{i}\triangle q_t\nabla_{\bar{i}}q_t
+\nabla_{\bar{i}}\triangle q_t\nabla_{i}q_t+R_{i\bar{j}}\nabla_{\bar{i}}q_t\nabla_j{q_t}.\nonumber
\end{eqnarray}
Hence  for any $p\geq 2$, it follows
\begin{eqnarray}\label{estimate-4.7}
\begin{aligned}
&\frac{4(p-1)}{p^2}\int_M|\nabla w_t^{p/2}|^2\omega_{g_t}^n=-\int_Mw_t^{p-1}\triangle q_t\omega_{g_t}^n\\
&=-\int_Mw_t^{p-1}(|\nabla\nabla q_t|^2+|\nabla\bar{\nabla}q_t|^2)\omega_{g_t}^n\\
&-2\text{Re} \int_M q_t^{p-1}\nabla_i\triangle q_t\nabla_{\bar{i}}q_t\omega_{g_t}^n
-\int_M q_t^{p-1}R_{i\bar{j}}\nabla_{\bar{i}}q_t\nabla_{j}q_t\omega_{g_t}^n.
\end{aligned}
\end{eqnarray}

On the other hand,  by Lemma  \ref{lem:perelman-1} and Theorem \ref{theorem-appen1},  we estimate
\begin{eqnarray*}
\begin{aligned}
&\quad -2\text{Re}\int_M w_t^{p-1}\nabla_i\triangle q_t\nabla_{\bar{i}}w_t\omega_{g_t}^n\\
&=-\text{Re}\int_Mw_t^{q-1}\nabla_i(|\nabla f_t|^2-2f+\triangle h_t)\nabla_{\bar{i}}q_t\omega_{g_t}^n\\
&=-\text{Re}\int_M(|\nabla f_t|^2-2f_t+\triangle h_t)(\frac{2(p-1)}{p}w_t^{\frac{p}{2}-1}\nabla_iw_t^{p/2}\nabla_{\bar{i}}q_t+w_t^{q-1}\triangle q_t)\omega_{g_t}^n\\
&\leq C(g)[\int_M \frac{2(p-1)}{p }  w_t^{\frac{p-1}{2}}|\nabla w_t^{p/2}|\omega_{g_t}^n+\int_M w_t^{p-1}|\triangle q_t|]\omega_{g_t}^n\\
&\leq\frac{p-1}{p^2} \int_M|\nabla w^{p/2}|^2\omega_{g_t}^n +   C(g)' p \int_M
w^{p-1}\omega_{g_t}^n
\end{aligned}
\end{eqnarray*}
and
\begin{eqnarray*}
\begin{aligned}
&\quad-\int_Mw_t^{p-1}R_{i\bar{j}}\nabla_{\bar{i}}q_t\nabla_{j}q_t\omega_{g_t}^n\\
&=-\int_Mw_t^p \omega_{g_t}^n-\int_Mw_t^{p-1}\nabla_i\nabla_{\bar{j}}h_t\nabla_{\bar{i}}q_t\nabla_{j}q_t\omega_{g_t}^n\\
&=-\int_Mw_t^p \omega_{g_t}^n+\int_Mw_t^{p-1}\nabla_{\bar{j}}h_t(\nabla_{\bar{i}}q_t\nabla_i\nabla_{j}q_t\omega_{g_t}^n+\triangle q_t\nabla_jq_t)\omega_{g_t}^n\\
&\quad+\frac{2(p-1)}{p}\int_Mw_t^{\frac{p}{2}-1}\nabla_iw_t^{p/2}\nabla_{\bar{j}}h_t\nabla_{\bar{i}}q_t\nabla_jq\omega_{g_t}^n\\
&\leq\int_Mw_t^{p-1}(|\nabla\nabla q_t|^2+\frac{1}{2}|\nabla\bar{\nabla}q_t|^2)\omega_{g_t}^n
+\frac{p-1}{p^2} \int_M|\nabla w_t^{p/2}|^2\omega_{g_t}^n\\
&+C(g)(p-1)\int_Mw_t^p\omega_{g_t}^n.
\end{aligned}
\end{eqnarray*}
Then  substituting  the above two inequalities into  (\ref{estimate-4.7}),   we get
\begin{eqnarray}
\int_M|\nabla w_t^{p/2}|^2\omega_{g_t}^n\leq C(g)(p-1)^2\int_Mw_t^{p-1}\omega_{g_t}^n,\hspace{0.2cm}\forall ~p\geq 2.\nonumber
\end{eqnarray}
By using Zhang's  Sobolev inequality [Zha], we deduce
\begin{eqnarray}\label{estimate-4.8}
\big(\int_M w_t^{p\nu}\big)^{1/\nu}\omega_{g_t}^n\leq C(g)C_s(q-1)^2\int_Mw_t^{p-1}\omega_{g_t}^n,\quad \forall ~p\geq 2,
\end{eqnarray}
where $\nu=\frac{n}{n-1}$.  To run the iteration we put $p_0=1$ and $p_{k+1}=p_k\nu+\nu$, $k\geq 0$.  Hence
\begin{eqnarray}
\|w_t\|_{L^{p_{k+1}}}&\leq&(CC_s)^{\frac{1}{p_k+1}}p_k^{\frac{2}{p_k+1}}\|w\|_{L^{p_k}}^{\frac{p_k}{p_k+1}}\nonumber\\
&\leq&(CC_s)^{\sum_{i=0}^{i=k}\frac{\nu^{k-i}}{p_k+1}}\prod_{i=0}^{i=k} p_i^{\frac{2\nu^{k-i}}{p_k+1}}\|w_t\|_{L^1}^{\prod\frac{p_i}{p_i+1}}\nonumber\\
&\leq&C(n,g)C_s^{\frac{n}{2}}\|w_t\|_{L^1}^{\gamma(n)}\nonumber
\end{eqnarray}
for a constant $\gamma(n)$ depending only on $n$, where we have used the fact $p_{k}\leq 2\nu^{k}$ for $k\geq 1$.  Therefore by (\ref{estimate-4.6}), we prove
$$\|w_t\|_{C^0}\leq C(n,g)C_s^{\frac{n}{2}}\|w_t\|_{L^1}^{\gamma(n)} \to 0,  ~\text{as}~t\to\infty.$$
\end{proof}

\vskip6mm
\section{Another version of the invariance $N_X(\omega_g)$  }

Let  $Y\in \eta_r(M)$ so that $\text{Im}(Y)$ generates an
one-parameter compact subgroup of $K$. Denote  $\mathcal{K}_{Y}$  to be a
class of $K_Y$-invariant  K\"{a}hler metrics in $2\pi c_{1}(M)$.
Then according to the proof of Proposition \ref{upperbound}, we
actually prove
\begin{eqnarray}
\begin{aligned}
\sup_{g\in\mathcal{K}_{Y}}\lambda(g)\le(2\pi)^{-n}[nV- \tilde F_Y(Y)-N_Y(c_1(M))].
\end{aligned}
\end{eqnarray}
Note that
$$\tilde F_{Y}(Y)=\int_{M}Y(h_g-\theta_{Y,\omega})e^{\theta_{Y,\omega_g}}\omega_g^{n}$$
and
 $$N_Y(c_1(M))= \int_{M}\theta_{Y,\omega_g}e^{\theta_{Y,\omega_g}}\omega_{g}^{n}
$$
 are both holomorphic invariances of $M$.  In this section, we want to show

\begin{proposition}\label{proposition-3} Let $H(Y)=\tilde F_Y(Y)+N_Y(c_1(M))$.
Then
$$\sup_{Y\in\eta_r(M)} H(Y)=N_X(c_1(M)),$$
where $X$ is  the extremal vector field as in Section 1. 
 \end{proposition}

\begin{proof}
Choose a constant $c_Y$ so that
$\hat\theta_{Y,\omega_g}=\theta_{Y,\omega_g}+c_Y$ satisfies a
normalization condition
\begin{eqnarray}
\begin{aligned}\label{normli-3}
\int_{M}\hat\theta_{Y,\omega_g}e^{h_g}\omega_g^{n}=0.
\end{aligned}
\end{eqnarray}
Then $\hat\theta_{Y,\omega_g}$ satisfies an equation
\begin{equation*}
\Delta \hat\theta_{X,\omega_g}+X(h_g)+\hat\theta_{X,\omega_g}=0
\end{equation*}
Thus using the integration by part, we have
\begin{equation*}\label{identity}
\tilde F_{Y}(Y)+\int_{M}\hat\theta_{Y,\omega_g}e^{\theta_{Y,\omega_g}}\omega_g^{n}=0
\end{equation*}
It follows
\begin{equation}\label{identity}
H(Y)=-c_YV=\int_{M}\theta_{Y,\omega_g}e^{h_g}\omega_g^{n}.
\end{equation}

We compute the first variation of $H(Y)$ in $\eta_r(M)$. By the
definition of $\theta_{Y+tY'}$,  we see that there exist constants
$b(t)$ such that $\theta_{Y+tY'}=\theta_{Y}+t\theta_{Y'}+b(t)$.
Since $\int_{M}e^{\theta_{Y+tY'}}\omega_g^{n}=V$, we have
$$e^{-b(t)}=\frac{1}{V}\int_{M}e^{\theta_{Y}+t\theta_{Y'}}\omega_g^{n}.$$
Thus we get
\begin{equation}
\frac{dH(Y+tY')}{dt}|_{t=0}=\int_{M}\theta_{Y'}e^{h_g}\omega_g^{n}-\int_{M}\theta_{Y'}e^{\theta_{Y}}\omega_g^{n}=\tilde F_{Y}(Y').
\end{equation}
 Therefore,   by [TZhu2], we see that there exists a unique critical  $X\in\eta_r(M)$ of
$H(\cdot)$ such that
 \begin{equation}\label{critical}
\tilde F_{X}(Y')= F_{X}(Y')\equiv 0,~~\forall~Y'\in\eta_r(M).
\end{equation}

Similarly, we have
$$\theta_{tY+(1-t)Y'}=t\theta_{Y}+(1-t)\theta_{Y'}+b(t)', ~\forall ~t\in [0, 1]$$
 for some constants $b(t)'$. Then
   \begin{eqnarray*}
\begin{aligned}
V=\int_{M}e^{\theta_{tY+(1-t)Y'}}\omega_g^{n}&=e^{b(t)'}\int_{M}e^{t\theta_{Y}+(1-t)\theta_{Y'}}\omega_g^{n}\\
&\leq e^{b(t)'}[t\int_{M}e^{\theta_{Y}}\omega_g^{n}+(1-t)\int_{M}e^{\theta_{Y'}}\omega_g^{n}]\\
&=e^{b(t)'}V.
\end{aligned}
\end{eqnarray*}
Thus  $b(t)'\geq 0$.  Consequently
 $$H(tX+(1-t)Y)\geq
tH(X)+(1-t)H(Y).$$
 This means that  $H(\cdot)$ is a concave
functional on $\eta_{r}(M)$. It follows that $X$ is a global
maximizer of $H(\cdot)$.  Therefore   we prove the proposition by using the fact $H(X)=N_X(c_1(M))$.

\end{proof}

\begin{corollary}\label{coro-2}  Let $\mathcal{K}_{K}$ be a class of
$K$-invariant  K\"{a}hler metrics in $2\pi c_{1}(M)$.  Suppose that
\begin{eqnarray}
\begin{aligned}
\sup_{g\in\mathcal{K}_{K}}\lambda(g)<
\inf_{Y\in\eta_{r}(M)}(2\pi)^{-n}[nV-F_Y(Y)-N_Y(c_1(M))].
\end{aligned}
\end{eqnarray}
 Then $(M,J)$ could not admit any  K\"ahler-Ricci soliton. Furthermore,    the modified Mabuchi's K-energy  could not be bounded from below.

\end{corollary}

\begin{proof}  The first part of corollary follows from Proposition \ref{proposition-3}  and Corollary \ref{coro-1}.  The second part  follows from Proposition \ref{proposition-3} and Theorem \ref{theorem-1}.

\end{proof}

The above corollary gives  a new obstruction to the existence of K\"ahler-Ricci solitons.

\vskip6mm
\section{ Proof of Theorem \ref{theorem-2}}

In this section,  we will modify the proof  of Main Theorem in [TZhu5] to prove Theorem \ref{theorem-2}.  The  proof in [TZhu5]  depends on a generalized   uniqueness theorem  for  K\"ahler-Einsteins recently  proved  by Chen and Sun  in [CS].
 Here we  avoid to use  their theorem
so that we can generalize the proof  to the case of K\"ahler-Ricci solitons by applying  Proposition \ref{prop-2}.
As in [TZhu5], we write an initial K\"ahler form $\omega_g$  of K\"ahler-Ricci  flow (2.1) by
$$\omega_g=\omega_\varphi=\omega_{g_{KS}}+\sqrt{-1}\partial\overline\partial\varphi\in 2\pi c_1(M)$$
for a K\"ahler potential $\varphi$ on $M$. We define a path of
K\"ahler forms
$$\omega_{g^s}=\omega_{g_{KS}}+s\sqrt{-1}\partial\overline\partial\varphi$$
and set
 \begin{eqnarray*}
\begin{aligned} I=\{ s \in ~[0,1]|~&(g_t^s; g^s)~\text{converges to a K\"ahler-Ricci soliton in }~C^\infty ~~\\
&\text{in sense of K\"ahler potentials}\}.
\end{aligned}
 \end{eqnarray*}
Clearly,  $I$ is not empty by the assumption of existence of K\"ahler-Ricci solitons on $M$.
We want to show that $I$ is in fact  both open and closed.  Then it follows that $I=[0,1]$. This  will finish the proof Theorem \ref{theorem-2}.

The openness of  $I$ is related to  the following   stability theorem of  K\"ahler-Ricc flow, which was proved in [Zhu2].

\begin{lemma}\label{stability} Let $(M,J)$ be a compact K\"ahler manifold
which admits a  K\"ahler-Ricci soliton $(g_{KS},X)$.    Let $\psi$ be a
K\"ahler potential of a $K_X$-invariant initial metric $g$ of (2.1).   Then there exists a small $\epsilon$
such that if
$$\|\psi\|_{C^{3}}\le \epsilon,$$
the solution $g(t,\cdot)$ of (2.1) will converge to a  K\"ahler-Ricci soliton with  respect to $X$ in
 $C^\infty$ in the sense of K\"ahler potentials.   Moreover, the convergence can be made
exponentially.

\end{lemma}

\begin{remark}\label{remark-3} Lemma \ref{stability} is still true if the $K_X$-invariant condition is  removed for the  initial metric $g$ (cf. [Zhu2]).  But we do not know  whether  the  convergence is
exponentially fast or not.
\end{remark}

\begin{proof}[Proof of openness of $I$]  Suppose that  $s_0\in I$.   Then  by the uniqueness of K\"ahler-Ricci solitons [TZhu1],  the flow $(g_t^{s_0}; \omega_{s_0})$   converges to $g_{KS}$ after a holomorphism
transformation in  $\text{Aut}_r(M)$.  Namely,  there exists a $\sigma\in   \text{Aut}_r(M)$ such that
$\sigma^\star\omega_{g_t^{s_0}}=\omega_{KS}+\sqrt{-1}\partial\overline\partial (\varphi_t^{s_0})_\sigma$ with property
$$\|(\varphi_t^{s_0})_\sigma\|_{C^k}\le C_k e^{-\alpha_k t},$$
where $C_k, \alpha_k>0$ are two uniform constants.
Then
 we can choose $T$ sufficiently large such that
$$\| (\varphi_t^{s_0})_\sigma\|_{C^{3}(M)}< \frac{\delta}{2},$$
where $\delta$ is a small number determined in Lemma  \ref{stability}.
Since the K\"ahler-Ricci flow is stable for any fixed finite time, there is a small $\epsilon >0$ such that
$$\| \varphi_T^{s} - (\varphi_T^{s_0})_\sigma\|_{C^{3}(M)}\,<\, \frac{\delta}{2}, ~\forall ~s\in [s_0,s_0+\epsilon],$$
where $\varphi_T^{s}$  is a K\"ahler potential of evolved K\"ahler metric $g_T^s$ of K\"ahler-Ricci flow $(g_t^s; \sigma^\star\omega_s)$ at time $T$.  Hence, we have
\begin{align}\| \varphi_T^{s}\|_{C^3(M)}\,<\, \delta,  ~\forall ~s\in [s_0,s_0+\epsilon].\end{align}
Then the flow $(g_t;  g_T^s)$ with initial $g_T^s$ will
converge to  a K\"ahler-Ricci soliton  in  $C^\infty$ according to Lemma \ref{stability}. This shows $s\in I$  for any  $s\in [s_0,s_0+\epsilon]$

\end{proof}

Let  $\varphi_t^{s}$  be  a family of  K\"ahler potentials of evolved K\"ahler metric $g_t^s$ of K\"ahler-Ricci flow $(g_t^s; \omega_s)$.   To make potentials  $\varphi^s_t$ more smaller to control,
 we need the following lemma, which was proved in [TZhu1].

\begin{lemma}\label{lemma-5.2} Let $M$ be a compact K\"ahler manifold
which admits a  K\"ahler-Ricci soliton $(g_{KS},X)$.    Let $\varphi$ be  a $K_X$-invariant
K\"ahler potential.   Then there exists a unique   holomorphism transformation $\sigma\in \text{Aut}_r(M)$
such that  $\varphi_\sigma\in \Lambda^{\perp}(\omega_{KS})$ with property
$$J(\varphi_\sigma)=\inf_{\tau\in\text{Aut}_r(M)}J(\varphi_\tau),$$
 where $\Lambda^{\perp}(\omega_{KS})$ is an orthogonal space to kernel space of linear operator  $(\Delta_{g_{KS}}+X+\text{Id})(\psi)$ and
 $$J(\varphi)=-\int_M \varphi e^{\theta_{X,\omega_{\varphi}}}\omega_{\varphi}^n+ \int_0^1\int_M \varphi e^{\theta_{X,\omega_{\lambda\varphi}}}\omega_{\lambda\varphi}^n\wedge d\lambda \ge 0.$$
   Moreover,
$$\|\sigma-\text{Id}\|\le C(\|\varphi\|_{C^5}),$$
where $\|\sigma-\text{Id}\|$ denotes the distance norm in Lie group $\text{Aut}_r(M)$ .

\end{lemma}

\begin{proof}[Proof of closedness  of $I$] By the openness of $I$,    we see that there exists a $\tau_0\le 1$ with
   $[0,\tau_0)\subset I$.  We need to show that $\tau_0\in I$.
In fact we want to prove that for any
$\delta>0$ there exists a large $T$ such that
\begin{align}\label{claim-5}\|(\phi^s_t)_{\sigma_{s,t}}\|_{C^{5}}\le \delta,~\forall~
t\ge T~~\text{and}~s< \tau_0,\end{align}
 where $\sigma_{s,t}$  are some holomorphisms in $\text{Aut}_r(M)$. We  will use an  argument by
contradiction as in  [TZhu5].
On contrary,
we can find a sequence of  evolved K\"ahler metrics
$g_{t_i}^{s_i}$ of K\"ahler-Ricci flows $(g_t^{s_i}; g^{s_i})$,  where $s_i\to
\tau_0$ and $t_i\to\infty$,  and a sequence of unique  holomorphisms $\sigma_{s_i,t_i}\in \text{Aut}_r(M)$ for pairs $(s_i,t_i)$
 such that $ (\phi^{s_i}_{t_i})_{\sigma_{s_i,t_i}}\in \Lambda^{\perp}(\omega_{KS})$ and
\begin{align} \| (\phi^{s_i}_{t_i})_{\sigma_{s_i,t_i}} \|_{C^{5}}\ge  \delta_0>0,\end{align}
for some constant $\delta_0$. Since the  K\"ahler-Ricci   flow $(g_t^{s_i}; g^{s_i})$ converges to some K\"ahler-Ricci soliton,   by the uniqueness of K\"ahler-Ricci solitons  [TZhu1],  the flow after a holomorphism
transformation in  $\text{Aut}_r(M)$   converges to $g_{KS}$.    So  we may further
assume that $ \phi^{s_i}_{t_i}$ satisfy
\begin{align} \|(\phi^{s_i}_{t_i})_{\sigma_{s_i,t_i}}   \|_{C^{5}}\le  2\delta_0.
\end{align}
Then there exists  a subsequence $(\phi^{s_i}_{t_i})_{\sigma_{s_i,t_i}}$ (still used by  $(\phi^{s_i}_{t_i})_{\sigma_{s_i,t_i}}$ )
of   $(\phi^{s_i}_{t_i})_{\sigma_{s_i,t_i}}$   converging to a potential $\phi_\infty\in \Lambda^{\perp}(\omega_{KS})$
with property
\begin{align}\label{contradiction-2} 2\delta_0\ge \|\phi_\infty\|_{C^{5}}\ge
\delta_0.\end{align}

We want to show that
 \begin{align}\label{identity-3}\lambda(\omega_{\phi_\infty)}=\lambda(g_{KS})=(2\pi)^{-n}(nV-N_X(c_1(M)).\end{align}
First we note that the modified $K$-energy is bounded from below since $M$ admits a  K\"ahler-Ricci soliton [TZhu2]. Then by  Proposition \ref{prop-2} and the monotonicity of $\lambda(g_t^{\tau_0})$,  we see  that
for any $\epsilon>0$, there exists a large  $T> 0$ such that
$$\lambda(g^{\tau_0}_t)\ge (2\pi)^{-n}(nV-N_X(c_1(M))-\frac{\epsilon}{2},~~\forall~t\ge T.$$
Since K\"ahler-Ricci flow is stable in finite time and $\lambda(g_t^s)$ is monotonic in $t$,
 there is a small $\delta> 0$ such that
for any $s\ge \tau_0-\delta$, we have
\begin{align}\lambda(g_t^s)\ge (2\pi)^{-n}(nV-N_X(c_1(M))-\epsilon,~~\forall~t\ge T.\end{align}
Since $s_i\to \tau_0$ and $t_i\to \infty$, we conclude that
$$\lim_{s_i\to \tau_0, t_i\to \infty} \lambda(\sigma_{s_i,t_i}^\star  g^{s_i}_{t_i})= \lim_{s_i\to \tau_0, t_i\to \infty}\lambda(g^{s_i}_{t_i})=(2\pi)^{-n}(nV-N_X(c_1(M)).$$
By the continuity of $\lambda(\cdot)$,  we will get  $(\ref{identity-3})$.

Now  by  Corollary 1.4 together with $(\ref{identity-3})$,  we see that $\omega_{\phi_\infty}$ is a  global maximizer of $\lambda(\cdot)$  in  $ \mathcal{K}_{X}$, so it is a  critical point of $\lambda(\cdot)$.
Then it is easy to show that $\omega_{\phi_\infty}$  a K\"ahler-Ricci soliton with respect to $X$ by computing the first variation of $\lambda(\cdot)$ as done in (\ref{variation}). 
Thus  by the uniqueness result for K\"ahler-Ricci solitons in [TZhu1], we get
$$\omega_{\phi_\infty}=\sigma^\star\omega_{KS},$$
 where $\sigma\in \text{Aut}_r(M)$.  Since $\phi_\infty\in \Lambda^{\perp}(\omega_{KS})$,  by Lemma \ref{lemma-5.2},  $\phi_\infty$ must be zero. This is a contradiction to
 $(\ref{contradiction-2})$.  The contradiction implies that (\ref{claim-5})  is true.

 By $(\ref{claim-5})$,  we see that  for any
$\delta>0$ there exists a large $T_0$  and $\sigma_0\in  \text{Aut}_r(M)$  such that
\begin{align}\label{claim-6}\|(\phi_{T_0}^{\tau_0})_{\sigma_{0}}\|_{C^{5}}\le  \delta.\end{align}
Then by Lemma \ref{stability}, the K\"ahler-Ricci  flow $(g_t; \omega_{(\phi_{T_0}^{\tau_0})_{\sigma_{0}} })$ converge to a K\"ahler-Ricci soliton.   Thus, we prove that  $\tau_0\in I$.

\end{proof}

\begin{remark}\label{remark-4} According to the proof of Theorem \ref {theorem-2} and Remark \ref{remark-3},   Theorem \ref {theorem-2}  will be  still true if the $K_X$-invariant condition is  removed for the  initial metric $g$ of
(2.1)  assuming that Conjecture (3.3) is true.
\end{remark}

\section{Appendix }

In [TZhu5],  it  was proved  the  minimizer  $f_t$ of $W(g_t, \cdot)$-functional associated to  evolved K\"ahler metric $g_t$   of  K\"ahler-Ricci  flow  (2.1) is  uniformly bounded (see also [TZha]).
In this appendix,  we show that the gradients of   $f_t$ are also uniformly bounded, and so are   $\triangle f_t$  by (\ref{EulerLagrange}).  Namely, we prove

\begin{theorem}\label{theorem-appen1}  There is a uniform constant $C$ such that
$$\|f_t\|+\|\nabla f_t\|+\|\triangle  f_t\|\le C, ~~\forall ~t>0.$$
\end{theorem}

We will derive  $\|\nabla f_t\|$ in Theorem \ref{theorem-appen1}  by  studying a  general  nonlinear elliptic equation as follows:
\begin{equation}\label{equation-appen1}
\triangle w(x)=w(x)F(x,w(x))
\end{equation}
where  the Laplace operator $\triangle $ is associated to a K\"ahler metric $g$ in $2\pi c_1(M)$ and  $F$ is a smooth function on $M\times\mathbb{R}^+$, which satisfies   a structure condition:
\begin{equation}\label{condition-appen1}
  -A-Bt^\alpha\leq F(\cdot,t)\leq H(t).
\end{equation}
Here  $0\le  A,B\le\infty  , 0\le \alpha<\frac{2}{n}$ are constants, and $H$ is a proper function on $\mathbb{R}^+$ which satisfies  a growth control at $0$:
\begin{equation}\label{e22.5}
\limsup_{t\rightarrow 0}\big(tH(t)\big)<\infty.
\end{equation}

\begin{lemma}\label{lemma-appen} Let $w$ is a positive solution of (\ref{equation-appen1}).  Then
\begin{equation}\label{estimate-appen1}
\|\nabla w\|_{C^0}\leq C(n)C_s^{\frac{n}{2}}\big(\|\nabla h\|_{C^0}+\|wF\|_{C^0}\big)^{n}\big(\int_M(1+|\nabla w|^2)dV_g\big)^{1/2},
\end{equation}
where $C_s$ is a Sobelev constant of $g$ and $h$ is a Ricci potential of $g$.
\end{lemma}

\begin{proof}  We will use the Moser's iteration to $L^p$-estimate of $|\nabla w|$.  By the Bochner formula, we have
\begin{eqnarray*}
\begin{aligned}
\triangle|\nabla w|^2&=|\nabla\nabla w|^2+|\nabla\bar{\nabla}w|^2+\nabla_i\triangle w\nabla_{\bar{i}}w+\nabla_iw\nabla_{\bar{i}}\triangle w+R_{i\bar{j}}\nabla_{\bar{i}}w\nabla_jw\\
&=|\nabla\nabla w|^2+|\nabla\bar{\nabla}w|^2+\nabla_i(wF)\nabla_{\bar{i}}w+\nabla_iw\nabla_{\bar{i}}(wF)+R_{i\bar{j}}\nabla_{\bar{i}}w\nabla_jw.
\end{aligned}
\end{eqnarray*}
Put $\eta=|\nabla w|^2+1$ . Then for $p\geq 2$, it follows
\begin{eqnarray}\label{estimate-appen4}
\begin{aligned}
&\frac{4(p-1)}{p^2}\int_M|\nabla \eta^{p/2}|^2dV_g\\
&=-\int_M \eta^{p-1}\triangle \eta dV_g \\
&=-\int_M \eta^{p-1}\big(|\nabla\nabla w|^2+|\nabla\bar{\nabla}w|^2\big)dV_g-\int_M\eta^{p-1}R_{i\bar{j}}\nabla_{\bar{i}}w\nabla_j w dV_g\\
&-\int_M\eta^{p-1}\big(\nabla_i(wF)\nabla_{\bar{i}}w+\nabla_{\bar{i}}(wF)\nabla_iw\big)dV_g.
\end{aligned}
\end{eqnarray}
The last term on the right hand side can be estimate as follows.
\begin{equation*}
\begin{aligned}
&-\int_M\eta^{p-1}\big(\nabla_i(wF)\nabla_{\bar{i}}w+\nabla_{\bar{i}}(wF)\nabla_iw\big)dV_g\\
&=\int_MwF\big(\nabla_i\eta^{p-1}\nabla_{\bar{i}}w+\nabla_{\bar{i}}\eta^{p-1}\nabla_iw+2\eta^{p-1}\triangle w\big)dV_g\\
&=\frac{2(p-1)}{p}\int_MwF\eta^{\frac{p}{2}-1}\big(\nabla_i\eta^{p/2}\nabla_{\bar{i}}w+\nabla_{\bar{i}}\eta^{p/2}\nabla_iw\big)dV_g\\
&+2\int_MwF\eta^{p-1}\triangle wdV_g.\\
\end{aligned}
\end{equation*}
Then
\begin{equation}\label{estimate-appen6}
\begin{aligned}
&-\int_M\eta^{p-1}\big(\nabla_i(wF)\nabla_{\bar{i}}wdV_g+\nabla_{\bar{i}}(wF)\nabla_iw\big)dV_g\\
&\leq\frac{2(p-1)}{p^2}\int_M|\nabla\eta^{p/2}|^2 dV_g
+2(p-1)p\int_M(wF)^2\eta^{p-2}(\eta-1)dV_g\\
&+\int_M\eta^{p-1}\big(\frac{(\triangle w)^2}{2n}+2n(wF)^2\big)dV_g\\
&\leq\frac{2(p-1)}{p^2}\int_M|\nabla\eta^{p/2}|^2dV_g+\frac{1}{2n}\int_M\eta^{p-1}(\triangle w)^2 dV_g\\
&+2[p(p-1)+n]\|wF\|_{C^0}^2\int_M\eta^pdV_g.
\end{aligned}
\end{equation}
For the second term on the right hand side, we note
$$R_{i\overline j}=g_{i\overline j}+h_{i\overline j}.$$
Then
\begin{eqnarray*}
\begin{aligned}
&-\int_M\eta^{p-1}R_{i\bar{j}}\nabla_{\bar{i}}w\nabla_jwdV_g\\
 &=\int_M\eta^{p-1}\nabla_i\nabla_{\bar{j}}h\nabla_{\bar{i}}w\nabla_jwdV_g-\int_M\eta^pdV_g\\
&=\frac{2(p-1)}{p}\int_M\eta^{p/2-1}\nabla_{\bar{j}}h\nabla_i\eta^{p/2}\nabla_{\bar{i}}w\nabla_jwdV_g\\
&+\int_M\eta^{p-1}\nabla_{\bar{j}}h\big(\triangle w\nabla_jw+\nabla_{\bar{i}}w\nabla_i\nabla_jw\big)dV_g-\int_M\eta^pdV_g.
\end{aligned}
\end{eqnarray*}
Thus
\begin{equation}\label{estimate-appen7}
\begin{aligned}
&-\int_M\eta^{p-1}R_{i\bar{j}}\nabla_{\bar{i}}w\nabla_jwdV_g\\
&\leq\frac{p-1}{p^2}\int_M|\nabla \eta^{p/2}|^2dV_g+p(p-1)\|\nabla h\|_{C^0}^2\int_M\eta^{p-2}(\eta-1)^2dV_g\\
&+\frac{1}{2n}\int_M\eta^{p-1}(\triangle w)^2dV_g+\frac{n}{2}\|\nabla h\|_{C^0}^2\int_M\eta^{p-1}(\eta-1)dV_g\\
&+\frac{1}{2n}\int_M\eta^{p-1}|\nabla\nabla w|^2dV_g+\frac{n}{2}\|\nabla h\|_{C^0}^2\int_M\eta^{p-1}(\eta-1)dV_g\\
&\leq\frac{p-1}{p^2}\int_M|\nabla \eta^{p/2}|^2dV_g+\frac{1}{2n}\int_M\eta^{p-1}(\triangle w)^2dV_g\\
&+\frac{1}{2n}\int\eta^{p-1}|\nabla\nabla w|^2 dV_g+[p(p-1)+n]\|\nabla h\|_{C^0}^2\int_M\eta^pdV_g.
\end{aligned}
\end{equation}

Substituting  (\ref{estimate-appen6}) and    (\ref{estimate-appen7})  into (\ref{estimate-appen4}),  we get
\begin{eqnarray*}
\begin{aligned}
&\frac{p-1}{p^2}\int_M|\nabla \eta^{p/2}|^2dV_g\\
&\leq-\int_M \eta^{p-1}\big(|\nabla\nabla w|^2+|\nabla\bar{\nabla}w|^2\big) dV_g\nonumber\\
&+\frac{1}{n}\int_M\eta^{p-1}(\triangle w)^2dV_g+\frac{1}{2n}\int_M\eta^{p-1}|\nabla\nabla w|^2dV_g\nonumber\\
&+[2(p-1)p+n]\big(\|\nabla h\|_{C^0}^2+\|wF\|_{C^0}^2\big)\int_M\eta^pdV_g\nonumber\\
&\leq C(n)p^2\big(\|\nabla h\|_{C^0}^2+\|wF\|_{C^0}^2\big)\int_M\eta^pdV_g\nonumber.
\end{aligned}
\end{eqnarray*}
It follows
\begin{equation}\nonumber
\int_M|\nabla \eta^{p/2}|^2dV_g\leq C(n)p^3\big(\|\nabla h\|_{C^0}^2+\|wF\|_{C^0}^2\big)\int_M\eta^pdV_g,\hspace{0.3cm}\forall~ p\geq 2.
\end{equation}
Therefore, by iteration, we derive
\begin{equation}\nonumber
\sup\eta\leq C(n)D^{n/2}\big(\int_M\eta^2dV_g\big)^{1/2},
\end{equation}
where $D=C_s\big(\|\nabla h\|_{C_0}^2+\|wF\|_{C^0}^2\big)$. This implies (\ref{estimate-appen1}) .
\end{proof}

\begin{proposition}\label{proposition-appen}
\begin{equation}
\|\nabla w\|_{C^0}\leq C(\|w\|_{L^2}),
\end{equation}
where  the constant $C$ depends  only on $n, C_s,A,B,\alpha,H$,  $\text{Vol}(g)$, $\|\nabla h\|_{C^0}$ and $\|w\|_{L^2}$.
\end{proposition}

\begin{proof} First we note that by using the standard Moser's iteration to equation
$$\triangle w(x)\ge   -A-Bw^\alpha,$$
 it is easy to see 
 $$\sup w\le C(1+\|w\|_{L^2}^\gamma)$$
 for some constants $C$ and $\gamma$ which depend only on  $n,C_s,  A,B,  \alpha,H$  and $\rm{Vol}(g)$. On the other hand, by  (\ref{equation-appen1}),
we have
$$\int_M|\nabla w|^2dV_g=-\int_Mw\triangle wdV_g=-\int_MwFdV_g.$$
Then we see  that $\|\nabla w\|_{L^2}$ is bounded by $\|w\|_{L^2}$. Thus  the proposition follows from  Lemma  \ref{lemma-appen}.

\end{proof}
 
Since $v_t=e^{-\frac{f_t}{2}}$ satisfies  (\ref{EulerLagrange-2})  which is  a type of  equation (\ref{equation-appen1}),  then  by Perelman's estimates  (d)  in Lemma \ref{lem:perelman-1} 
and Zhang's estimate for Sobelev constants associated to $g_t$ in [Zha] together with $C^0$-estimate for $f_t$ in [TZhu5],   we obtain   a uniform  gradient estimate  for $v_t$    from
Proposition \ref{proposition-appen}, and so for $f_t$ . By  equation (\ref{EulerLagrange}), we also derive a  a uniform  Laplacian estimate for $f_t$. 
 Thus Theorem  \ref{theorem-appen1} is true.    Theorem  \ref{theorem-appen1}  will be used in Section 4.

\vskip6mm

\bibliographystyle{amsplain}

\begin{thebibliography}{10}
\bibitem {Pe}[Pe] G. Perelman, The entropy formula for the Ricci flow and its geometric applications, arxiv:math.DG/0211159.
\end{thebibliography}

\end{document}